\documentclass[11pt]{article}
\usepackage[ansinew]{inputenc}
\usepackage{latexsym}
\usepackage{textcomp}
\usepackage{color}
\usepackage{amscd}

\input{amssym.def}
\input{amssym}

     \newcommand{\fa}{\goth{a}}

     \newcommand{\fo}{\goth{o}}
     \newcommand{\fp}{\goth{p}}
     
     \newcommand{\fr}{\goth{r}}

     \newcommand{\fO}{\goth{O}}
     \newcommand{\fP}{\goth{P}}

     \newcommand{\Ll}{\Bbb{L}}  

     \newcommand{\Q}{\Bbb{Q}}
     
     \newcommand{\Z}{\Bbb{Z}}

    \newcommand{\ol}[1]{\overline{#1}}
    \newcommand{\ul}[1]{\underline{#1}}
    
    \newcommand{\ti}[1]{\tilde{#1}}
    \newcommand{\group}[1]{\langle{#1}\rangle}

    \newcommand{\ot}{\otimes}
    \newcommand{\me}{^{-1}}
    \newcommand{\mal}{^{\times}}
    
    \newcommand{\df}{\stackrel{\mathrm{def}}{=}}
    \newcommand{\implies}{\Longrightarrow}
    \newcommand{\mr}{\mathrm}
    \newcommand{\clo}{^{\mr{c}}}

    \newcommand{\zl}{{\Bbb{Z}_l}}
    
    \newcommand{\ql}{{\Bbb{Q}_l}}

    \newcommand{\into}{\rightarrowtail}
    \newcommand{\onto}{\twoheadrightarrow}
    \newcommand{\lto}{\longrightarrow}
    \newcommand{\da}{\downarrow}

    \newcommand{\pht}{\phantom}

    \def\daz#1{#1\da\pht{#1}}

    \newcommand{\sda}{\da_{\raisebox{-1.4mm}{$\hsp{-1}\check{}$}}}

    \newcommand{\ga}{\gamma}
    \newcommand{\Ga}{\Gamma}

    \newcommand{\ze}{\zeta}
    \newcommand{\la}{\lambda}
    \newcommand{\La}{\Lambda}

    \newcommand{\om}{\omega}
    \newcommand{\al}{\alpha}
    \newcommand{\ve}{\varepsilon}

    \newcommand{\noi}{\par\noindent}
    \newcommand{\sn}{\par\smallskip\noindent}
    \newcommand{\mn}{\par\medskip\noindent}
    \newcommand{\bn}{\par\bigskip\noindent}
    \newcommand{\bbn}{\par\bigskip\bigskip\noindent}
    
    \newcommand{\Section}[2]{\bbn {\large #1\,. \ {\sc #2}}
                             \nopagebreak
                             \nz}

    \newcommand{\nf}[2]{\\[1.5ex]
                        \bmp{1cm}
                         (#1)
                        \emp 
                        \bmp{13.5cm}
                         \bct
                          $#2$
                         \ect
                        \emp\\[1.5ex]
         }

    \newcommand{\sss}{\scriptstyle}
    \newcommand{\nz}{\\[1ex]}
    
    \newcommand{\vsp}[1]{\noi
           \vspace*{#1mm}
      \noi$\!\!$}
    \newcommand{\hsp}[1]{\hspace*{#1mm}}

    \newcommand{\mmargin}{
     \textheight 230truemm
     \textwidth 155truemm
     \topmargin -10truemm
     \oddsidemargin 5truemm
     \evensidemargin 5truemm
     }

    \newcommand{\bmp}{\begin{minipage}}
    \newcommand{\emp}{\end{minipage}}
    \newcommand{\btb}{\begin{tabular}}
    \newcommand{\etb}{\end{tabular}}
    \newcommand{\barr}{\begin{array}}
    \newcommand{\earr}{\end{array}}
    \newcommand{\bit}{\begin{itemize}}
    \newcommand{\eit}{\end{itemize}}
    \newcommand{\ben}{\begin{enumerate}}
    \newcommand{\een}{\end{enumerate}}
    \newcommand{\bct}{\begin{center}}
    \newcommand{\ect}{\end{center}}
    \newcommand{\bfr}{\begin{flushright}}
    \newcommand{\efr}{\end{flushright}}
    \newcommand{\bea}{\begin{eqnarray*}}
    \newcommand{\eea}{\end{eqnarray*}}
    \newcommand{\bqo}{\begin{quote}}
    \newcommand{\eqo}{\end{quote}}
    \newcommand{\bdc}{\begin{description}}
    \newcommand{\edc}{\end{description}}
    \newcommand{\bdia}{\begin{CD}}
    \newcommand{\edia}{\end{CD}}

    \definecolor{light}{gray}{.3}

    \newcommand{\Hom}{\mathrm{Hom}}

    \newcommand{\tr}{\mathrm{tr\,}}
    
    \newcommand{\coker}{\mathrm{coker\,}}
    \newcommand{\ind}{\mathrm{ind\,}}
    \newcommand{\infl}{\mathrm{infl\,}}
    \newcommand{\im}{\mathrm{im\,}}
    \newcommand{\mod}{\mathrm{\ mod\ }}
    \newcommand{\res}{\mathrm{res\,}}
    
    \newcommand{\Det}{\mathrm{Det\,}}
    \newcommand{\sr}[2]{{\,\stackrel{#1}{#2}\,}}
    \newcommand{\fra}[2]{{\,\frac{#1}{#2}\,}}

    \newcommand{\theorem}{\sn
           \bdc
           \item[{\sc Theorem.}] \em }

    \newcommand{\Stop}{\edc \sn\rm} 

    \newcommand{\Lemma}[1]{\sn
           \bdc
           \item[{\sc Lemma {#1}.}] \em }
    
    \newcommand{\Proposition}[1]{\sn
           \bdc
           \item[{\sc Proposition {#1}.}] \em }
    \newcommand{\corollary}{\sn
           \bdc
           \item[{\sc Corollary.}] \em }

    \newcommand{\proof}{{\sc Proof.} \ }
    \newcommand{\remark}{\sn{\sc Remark.} \ }
    
    \newcommand{\definition}{\sn
           \bdc
           \item[{\sc Definition.}] \em }
    
    \newcommand{\bmx}{\left(\barr}
    \newcommand{\emx}{\earr\right)}

\mmargin

\def\frs{{\mr{Fr}_{\fp'}}}
\def\eulerk{{E(\fp,K/k)}}
\def\eulers{{E(\fp',K/k')}}
\def\euler1{{E(\fp_1,K/k_1)}}
\def\eulert{{E(\fp_t,K/k_t)}}
\def\eulerts{{E(\fp_t',K/k_t')}}
\def\cT{{{\cal{T}'}}}

\def\gs{{G'}}
\def\gi{{G}}
\def\lgi{{\La \gi}}
\def\lbgi{{\La_\bullet\gi}}
\def\lwgi{{\La_\wedge\gi}}
\def\gak{{\Ga_k}}

\def\lcgak{{\La\clo\gak}}
\def\lcwgak{{\La\clo_\wedge\gak}}
\def\lbgak{{\La_\bullet\gak}}
\def\lcbgak{{\La\clo_\bullet\gak}}
\def\lwgak{{\La_\wedge\gak}}
\def\lga{{\La\Ga}}
\def\lbga{{\La_\bullet\Ga}}
\def\lwga{{\La_\wedge\Ga}}
\def\lcga{{\La\clo\Ga}}
\def\lwgab{{\La_\wedge G^\mr{ab}}}

\def\lbgs{{\La_\bullet G'}}
\def\qgi{{{\cal{Q}}\gi}}
\def\qga{{{\cal{Q}}\Ga}}
\def\qgak{{{\cal{Q}}\gak}}

\def\qcgak{{{\cal{Q}}\clo\gak}}
\def\qcwgak{{{\cal{Q}}\clo_\wedge\gak}}
\def\qwgak{{{\cal{Q}}_\wedge\gak}}
\def\qwgi{{{\cal{Q}}_\wedge\gi}}
\def\qwga{{{\cal{Q}}_\wedge\Ga}}

\def\qwtg{{{\cal{Q}}_\wedge\ti G}}
\def\qd1{{{\cal{Q}}D_1}}
\def\qdi1{{{\cal{Q}}(D_1/I_1)}}
\def\qgs{{{\cal{Q}}\gs}}
\def\HOM{{\mr{HOM}}}
\def\homs{{\mr{Hom}^\ast}}
\def\rlgi{{R_l(\gi)}}

\def\giab{{G^\mr{ab}}}
\newcommand{\LL}{{{\mbox{{\boldmath$\mr{L}$}}}}} 
\def\tr{{\mr{Tr}}}

\def\Det{{\mr{Det}}}
\def\defl{{\mr{defl}}}
\def\deflt{{\mr{defl}_{G}^{\ti G}}}
\def\deflab{{\defl_{G}^{G^\mr{ab}}}}
\def\lki{{L_{K/k}}}

\def\lks{{L_{K/k'}}}

\def\Ki{{K}}
\def\ki{{k_\infty}}
\def\lamki{{\la_{\Ki/k}}}

\def\kab{{K_\mr{ab}}}

\def\tki{{t_{\Ki/k}}}
\def\tkt{{t_{\ti K/k}}}
\def\tkab{{t_{\kab/k}}}
\def\tkiga{{(\tki|\Ga)}}

\def\tkaga{{(\tkab|\ti\Ga)}}
\def\plc{{\psi_l\chi}}
\def\plcs{{\psi_l\chi'}}
\def\plcj{{\psi_l\chi_j}}
\def\ver{{\mr{ver}}}
\def\cver{{\chi'\circ\ver}}
\def\csver{{\chi'\circ\ver}}
\def\tga{{(t|\Ga)}}
\def\mchi{{\fra1l\log\fra{\lki(\chi)^l}{\Psi\lki(\plc)}}}
\def\giga{{\fra1{[\gi:\Ga]}}}

\def\indpsi{{\sum_{i=0}^{l-1}(\chi'\circ\ver)\om^i-l(\chi'\circ\ver)}}
\def\chisa{{{\chi'}^{\hat A}}}
\def\frp{{\mr{Fr}_\fp}}

\def\fr1{{\mr{Fr}_{\fP_1}}}
\def\np1{{\mr{N}{\fp_1}}}

\def\lkis1{{L_{\Ki/k,S_1}}}
\def\lkikss1{{L_{\Ki/k',S_1'}}}

\begin{document}
\title{{Non-abelian pseudomeasures and congruences between\\ abelian Iwasawa $L$-functions}}
 \author{Jürgen Ritter \ $\cdot$ \ Alfred Weiss \
 \thanks{We acknowledge financial support provided by NSERC and the University of Augsburg.}
 }
\date{\pht\today}

\maketitle \bct\footnotesize{Dedicated to Professor J.-P.~Serre on
his 80$^\mr{th}$ birthday}\ect
 \medskip
 \bct
 \bmp{12,5cm}
 {\small
 \noi
 {\sc Abstract\,.} \  The paper starts out from pseudomeasures (in
 the sense of Serre) which hold the arithmetic properties of the
 abelian $l$-adic Artin $L$-functions over totally real number fields.
 In order to generalize to non-abelian $l$-adic $L$-functions,
 these abelian pseudomeasures must satisfy congruences which are
 introduced but not yet known to be true. The relation to the
 ``equivariant main conjecture'' of Iwasawa theory is discussed.

 }

 \emp
 \ect

 \bbn
 Fix an odd prime number $l$ and a finite field extension $k/\Q$
 with $k$ totally real. Let $\ki$ be the cyclotomic
 $\zl$-extension of $k$ and $\Ki\supset\ki$ be a totally real
 Galois extension of $k$ with Galois group $\gi$ and so that
 $[\Ki:\ki]$ is finite. Setting $H=G_{\Ki/\ki}$ and
 $\gak=G_{\ki/k}$, we get the group extension $1\to
 H\to\gi\to\gak\to1$\,.
 \mn Let first $\gi$ be abelian. We consider the group algebra
 $\qgi=\mr{Quot}(\lga)[H]$ which results from a splitting of the above
 group extension, with $\lga\ (\simeq\zl[[T]])$
 the Iwasawa algebra of a preimage $\Ga$ of
 $\gak$ in $\gi$. From Serre's interpretation [Se] of the work
 [DR] of Deligne and Ribet on abelian $L$-functions over totally
 real fields it follows that there is a unique element
 $\lamki\in\qgi$ that encodes all the $l$-adic $L$-functions
 $L_l(s,\chi)$ of $\ki/k$ for the characters $\chi$ of $\gi$ with
 open kernel. In fact, choosing a preimage $\ga\in\Ga$ of a
 generator $\ga_k\in\gak$ and identifying $\ga-1$ with the
 variable $T$ (in $\lga\simeq\zl[[T]])$, this {\em pseudomeasure}
 $\lamki$ admits a power series expansion
 $$\lamki=\sum_{m\ge-1}a_mT^m\ \mr{with}\ a_m\in\zl[H]$$ which has
 the property that
 $$\sum_{m\ge-1}\chi(a_m)(\chi(\ga)u^n-1)^m=L_l(1-n,\chi)\
 \mr{for}\ n\ge1\ \mr{and\ all}\ \chi\ ,$$ where $u\in1+l\zl$ is
 determined by the action of $\ga_k$ on the $l$-power roots of
 unity. In particular, $\lamki\in(\lbga)[H]$, with $_\bullet$
 denoting localization
 \footnote{i.e., inverting all elements in $\La\Ga\setminus l\cdot\La\Ga$} at the prime ideal
 $l\La\Ga$ of $\La\Ga$;
 moreover, $\lamki$ is a unit in
 $\qgi$. If Iwasawa's $\mu$-invariant of
 $\Ki/k$ vanishes, then
 \footnote{See [RW3, proof of Corollary to Theorem 9]. The $\mu$-invariant is that of
 the Galois group $X=G(M/K)$ of the maximal abelian $l$-extension $M$ of
 $K$ unramified outside $l$.}
 $\lamki$ is even a unit in $\lbgi$.
 \mn Though the pseudomeasure carries all the information about
 the abelian $l$-adic $L$-functions, its arithmetic properties are
 not well understood at present as far as we know.
 \mn We now drop the assumption that $\gi$ is abelian. Then the
 total ring $\qgi$ of fractions \footnote{we invert all central
 regular elements of $\lgi$} of the Iwasawa algebra
 $\lgi=\zl[[\gi]]$ of $\gi$ is a finite dimensional semisimple
 algebra over $\qga$, for every central open subgroup $\Ga$ of
 $\gi$ which is isomorphic to $\zl$. Moreover, there is a
 homomorphism $$\Det:K_1(\qgi)\to\Hom(\rlgi,(\qcgak)\mal)$$
 taking $x\in K_1(\qgi)$ to the map $\Det\,x$ which
 assigns an element in $(\qcgak)\mal\df(\ql\clo\ot_\ql\qgak)\mal$
 to each character $\chi$ of $\gi$ (with open kernel and with
 values in a fixed algebraic closure $\ql\clo$ of $\ql$). Above,
 $\rlgi$ is the $\Z$-span of the $\ql\clo$-irreducible characters
 $\chi$ with open kernel (see also [RW2, p.588]). Given a finite set $S$ of
 primes
 of $k$ which contains all primes above $\infty$ and $l$, Greenberg [Gr]
 has generalized $l$-adic $L$-functions $L_{l,S}(s,\chi)$ to non-abelian characters
 $\chi$ of $\gi$ (with open kernel), and these in turn define the
 Iwasawa $L$-function $\chi\mapsto L_{\Ki/k}(\chi)$ (see [RW2,
 p.563]) which belongs to $\Hom(\rlgi,(\qcgak)\mal)$.
 \mn The natural question arises whether $L_{\Ki/k}$ has a
 preimage in $K_1(\qgi)$. Any such \footnote{which should be
 unique according to a conjecture of Suslin (see [RW2, p.565,
 Remark (E)])} may be regarded as a non-abelian analogue of Serre's
 abelian pseudomeasure. If $\gi$ is abelian,
 $K_1(\qgi)=(\qgi)\mal$ and we are back in Serre's situation. If
 Iwasawa's $\mu$-invariant of $\Ki/k$ vanishes, the abelian case
 hints at finding a $\la$ in $K_1(\lbgi)$
 \,\footnote{where $\La_\bullet G$ is obtained by inverting all central elements of $\La G$
 which are regular in $\La G/l\La G$}
 with $\Det\,\la=L_{\Ki/k}$;
 its image in $K_1(\qgi)$ is then a non-abelian pseudomeasure for $\Ki/k$\,.
 \mn The purpose of this paper is to formulate conditions which
 guarantee the existence of non-abelian pseudomeasures in the case when $\mu(\Ki/k)=0$. These are
 stated in terms of congruences between values of Iwasawa
 $L$-functions and lead to hypothetical new congruences between
 abelian pseudomeasures.
 \mn Here is a short outline of the contents of the paper. In \S1
 we review results from [RW2,3] and also provide some notation
 that is used throughout. In this and in the later sections we
 restrict ourselves to pro-$l$ groups $\gi$, though this would not
 be necessary. We also assume that Iwasawa's $\mu$-invariant of
 $\Ki/k$ vanishes (it then does so for all intermediate extensions
 $\ti K/k$ and $\Ki/k'$, with $\ti K$ and $k'$ the fixed fields of a finite normal
 subgroup $N\lhd\gi$ respectively an open subgroup $G'\le\gi$, see [RW3,
 footnote 1]). The next two sections introduce two kinds of congruences between Iwasawa
 $L$-functions, which hold precisely when $\lki\in \Det\, K_1(\lbgi)$\,. The second
 kind is reformulated in terms of congruences between abelian pseudomeasures. Section 4
 discusses the first kind for
 special pro-$l$ groups $\gi$ (those with an
 abelian subgroup of index $l$) and reduces them to the above mentioned
 congruences between abelian pseudomeasures if these groups $\gi$ have nilpotency class 2.
 \mn Finally, in an
 appendix, we prove that the existence of a pseudomeasure in
 $K_1(\qgi)$ does not depend on the size of $S$ as long as $S$
 contains all primes dividing $l$ or $\infty$. The
 appendix includes moreover a short review of the ``equivariant main
 conjecture'' of Iwasawa theory, as introduced in [RW2], and its
 equivalence to the existence of non-abelian pseudomeasures in $K_1(\qgi)$. Here, the set $S$ must
 be {\em sufficiently large}\,, i.e., contain all infinite primes and all
 primes of $k$ whose ramification index in $\Ki/k$
 is divisible by $l$ \footnote{so, in particular,
 $S$ contains the primes dividing $l$}. As a corollary, the ``equivariant main
 conjecture'' is independent of the choice of $S$.
 \bqo{\small{\bf Added in proof\,:} From the referee's report we
 have learned of the manuscript {\em Iwasawa theory of totally real
 fields for Galois extensions of Heisenberg type} by K.~Kato,
 which is closely related to this paper. We would like to thank
 Professor Kato for sending us the ``Very preliminary version''
 of January 15, 2007.
 \sn In this manuscript, Kato sketches a proof of the main
 conjecture of its title, which (in Lie dimension 1) is equivalent to the ``equivariant
 main conjecture'' (see \S5), in the special case of Galois groups
 of Heisenberg type. Very roughly, the idea is to identify certain
 congruences between abelian pseudomeasures (we would say), and
 then to verify them by using methods from [DR].
 \sn When the Heisenberg situation overlaps with \S4, these
 certain congruences are a variation of the single hypothetical
 congruence \ $\la_{K/k'}\equiv\ver(\la_{K_\mr{ab}/k})\mod\cT$ \ which
 is equivalent to the main conjecture by Propositions 3.2 and 4.4.
 So long as the relevant congruences can be explicitly stated,
 there is cause for optimism that Kato's ideas will permit a proof
 of the main conjecture. In this connection we would like to
 mention that the nilpotency class 2 assumption in Proposition 4.4
 has meanwhile been eliminated.}\eqo
 \Section{1}{Notation and basic results}
 Whenever we write $\Ga$ we mean an open subgroup of $\gi$ which
 is isomorphic to $\zl$; whenever we say ``character'' we mean a
 $\ql\clo$-valued character of $\gi$ with open kernel. Such
 characters are said to be of type W, and always denoted by
 $\rho$, if $\rho$ is $\ql\clo$-irreducible and satisfies
 $\rho(h)=1$ for $h\in H$. Each $\rho$ determines the automorphism
 $\rho^\sharp$ of the field $\qcgak$, induced by
 $\rho^\sharp(\ga)=\rho(\ga)\ga$ for $\ga\in\gak$. The $l$-th
 Adams operation on $\rlgi$ is denoted by $\psi_l$, so
 $(\plc)(g)=\chi(g^l)$ for $g\in\gi$. And $\Psi$ is the endomorphism
 of $\lga$ induced by $\Psi(\ga)=\ga^l$ for $\ga\in\Ga$. Finally,
 set $\lcga=\zl\clo\ot_\zl\lga$, with $\zl\clo$ the ring of
 integers in $\ql\clo$.
 \mn With this notation we specify the subgroup
 $$\HOM(\rlgi,(\lcgak)\mal)\quad\mr{of}\quad\Hom(\rlgi,(\lcgak)\mal)$$ to
 consist of all homomorphisms $f:\rlgi\to(\lcgak)\mal$ satisfying
 \ben\item equivariance with respect to the natural action of
 $G_{\ql\clo/\ql}$ on $\rlgi$ and $\lcgak$\,,
 \item $f(\chi\ot\rho)=\rho^\sharp(f(\chi)) \quad(\forall\,\chi$ and all $\rho$ of type W)\,,
 \item $\fra{f(\chi)^l}{\Psi f(\plc)}\equiv 1\mod
 l\lcgak\quad(\forall\,\chi)$\ .\een
 {\sc Fact 1\,.}\quad{\em The determinant
 $\Det:K_1(\qgi)\to\Hom(\rlgi,(\qcgak)\mal)$ induces
 $$\Det:K_1(\lgi)\to\HOM(\rlgi,(\lcgak)\mal)\ .$$ Above,
 $\lgi,\lcgak$ may be replaced by $\lbgi,\lcbgak$ or by $\lwgi,\lcwgak$,
 where $\La_\wedge-$ is the $l$-adic completion of $\La_\bullet-$\,.}
 \mn NB \ The appearance of $\La_\bullet-$ has already been
 observed in the abelian case. Since we will have to work with
 logarithms, we need to complete in order to guarantee the
 convergence of the logarithm series.
 \mn For the actual definition of
 $\Det:K_1(\qgi)\to\Hom(\rlgi,(\qcgak)\mal)$ see [RW2, p.558]; for
 {\sc Fact 1} compare [RW3, Lemma 2, Propositions 4 and 11].
 Explicitly, we will use the formula $$(\Det\, g)(\chi)=\chi(g)\ol g
 \quad \mr{for}\quad g\in\gi \quad \mr{with\ image} \quad \ol g\in\gak$$
 whenever $\chi$ is irreducible of degree 1 [RW3, Proposition 5].
 \mn Introducing \ $T(\qgi)=\qgi/[\qgi,\qgi]$\,, with
 $[\qgi,\qgi]$ the additive subgroup generated by all
 $ab-ba\,,\,a,b\in\qgi$\,, we obtain the trace isomorphism
 \footnote{The superscript $\ast$ means HOM without condition 3. In [RW3]
 $\HOM,\,\La$ and functions $f$ often carry superscripts such as $\ast,N,\fO,\mr{Fr}$\,:
 since we are going to restrict $\gi$ to be a
 pro-$l$ group, these exponents should just be ignored.}
 $$\tr:T(\qgi)\to\homs(\rlgi,\qcgak)$$ given by $\tr(\tau
 g)(\chi)=\chi(g)\ol g$\,, where $\tau g$ denotes the image of
 $g\in\gi$ in $T(\qgi)$ (compare [RW3, Lemma 6, Proposition 3]).
 We remark that the natural map $T(\lgi)\to T(\qgi)$ is injective.
 \mn {\em From now on $\gi$ is always a pro-$l$ group.} Define
 $$\LL:\HOM(\rlgi,(\lcgak)\mal)\to\homs(\rlgi,\qcgak)\ \ \mr{by}\
 \  f\mapsto[\chi\mapsto\fra1l\log\fra{f(\chi)^l}{\Psi f(\plc)}]\ .$$
 {\sc Fact 2\,.}\quad {\em The homomorphism $\LL$ induces a unique
 homomorphism \ $\Ll:K_1(\lgi)\to T(\qgi)$ \ making the square
 $$\barr{ccc} K_1(\lgi)&\sr{\Ll}{\lto}&T(\qgi)\\
 \daz{\Det}&&\daz{\tr,\simeq}\\
 \HOM(\rlgi,(\lcgak)\mal)&\sr{\LL}{\lto}&\homs(\rlgi,\qcgak)\earr$$
 commute. This $\Ll$ takes $K_1(\lgi)$ into $T(\lgi)$. The same holds
 with $\lgi,\lcgak,\qgi,\qcgak$ replaced by
 $\lwgi,\lcwgak,\qwgi,\qcwgak$, respectively \footnote{Contrary to
 [RW3], we prefer to write $\La_\wedge-$ rather than
 $(\La-)_\wedge$ and the same with the total ring of fractions
 ${\cal{Q}}_\wedge-$ of $\La_\wedge-$.}.}
 \mn For this, see [RW3, Theorem 8, Proposition 11].
 \mn We finally turn to the Iwasawa $L$-function $\lki$. For this,
 we first fix a finite set $S$ of places of $k$ which contains all
 archimedean places and all places of $k$ above $l$. For $\chi\in\rlgi$, let
 $L_{l,S}(s,\chi)$ be the $l$-adic Artin
 $L$-function with respect to $S$, as defined in [Gr]. By [C-N] or [DR]
 there is an element
 $G_{\chi,S}(T)\in\ql\clo\ot_\ql\mr{Quot}\zl[[T]]$ such that \
 $L_{l,S}(1-s,\chi)=\fra{G_{\chi,S}(u^s-1)}{H_\chi(u^s-1)}$\,,  with
 the 1-unit $u\in\zl$ that has already appeared in the
 introduction and with $H_\chi(T)=\chi(\ga_k)(1+T)-1$ or $=1$
 according as $H\le\ker\chi$ or not.
 \mn With this notation $\lki$ is
 defined by \
 $\lki(\chi)=\frac{G_{\chi,S}(\ga_k-1)}{H_\chi(\ga_k-1)}$\,. Note
 that although $\lki$ depends on $S$ we suppress this dependence in our notation; note also
 that $\lki$ is independent of a special choice of $\ga_k$ [RW2,
 Proposition 11].
 \mn {\sc Fact 3\,.}\quad {\em If $S$ is sufficiently large and if Iwasawa's $\mu$-invariant of $\Ki/k$
 vanishes,  \ $\lki\in\HOM(\rlgi,(\lcbgak)\mal)\,.$ Moreover,
 $\lki\in\Det\,K_1(\lbgi)$ if, and only if, $\lki\in\Det\,K_1(\lwgi)$\,.}
 \mn For this see [RW3, Corollary of Theorem 9] and compare Theorems A
 in [RW3] and [RW4]. We remark
 that the assumption $\mu=0$ is independent of the size of $S$
 [NSW, (11.3.6)],
 because no prime of $k$ splits completely in the cyclotomic
 $\zl$-extension $\ki$ of $k$.
 \sn {\em From now on we assume $\mu=0$ for the extension
 $\Ki/k$.}
 \Section{2}{The logarithmic pseudomeasure and integrality}
 We start out from the diagram shown in {\sc Fact 2}, but
 now read in the completed situation\,:
 $$\barr{ccc} K_1(\lwgi)&\sr{\Ll}{\lto}&T(\qwgi)\\
 \daz{\Det}&&\daz{\tr,\simeq}\\
 \HOM(\rlgi,(\lcwgak)\mal)&\sr{\LL}{\lto}&\homs(\rlgi,\qcwgak)\,.\earr$$
 \definition $\tki\in T(\qwgi)$ is the unique element satisfying \
 $\tr(\tki)=\LL(\lki)$\,; we call it the logarithmic
 pseudomeasure of $K/k$.\Stop
 Recall that we are assuming $\mu=0$ for the pro-$l$ extension $\Ki/k$ and note that then
 $\tki=\Ll(\lamki)$ if $\gi$ is abelian (compare [RW3, Corollary to Theorem 9]).
 \Lemma{2.1} If $N$ is a finite normal subgroup of $\gi$ with fixed
 field $\ti K$, then $\deflt\tki=\tkt$\,, where $\ti G=G_{\ti K/k}\,(\simeq\gi/N)$.\Stop
 To see this, we recall that on the $K_1$-level deflation is
 induced by $\La_\wedge\ti G\ot_\lwgi-$\,, on the $T$-level by
 $\gi\onto\ti G$\,, and on the Hom-level by $f\mapsto[\ti\chi\mapsto
 f(\infl_{\ti G}^\gi\ti\chi)]$ for $\ti\chi\in R_l(\ti G)$\,. By
 [RW2, Lemma 9], Det and $\deflt$ commute and it is shown below that
 also Tr and $\deflt$ commute. Finally, $\LL$ and $\deflt$ commute because
 $\psi_l$ and $\infl_{\ti G}^\gi$ do. Hence the lemma follows from
 $\deflt(\lki)=L_{\ti K/k}$ (see [RW2, p.563]).
 \sn In order to check commutativity of the diagram
 $$\barr{ccc}T(\qwgi)&\sr{\defl}{\lto}&T(\qwtg)\\
 \daz{\tr}&&\daz{\ti\tr}\\
 \homs(\rlgi,\qcwgak)&\sr{\defl}{\lto}&\homs(R_l(\ti
 G),\qcwgak)\,,\earr$$
 we choose a central open subgroup $\Ga$ in $\gi$ and representatives
 $g_i\in\gi,\,1\le i\le s\,,$ of the conjugacy classes of the
 finite group
 $\gi/\Ga$. Writing $\ti{\pht{x}}$ for the map $\gi\to\ti G$, take
 $x\in T(\qwgi)$ and a character $\ti\chi$ of $\ti
 G/\ti\Ga\simeq\gi/\Ga N$. Writing $x=\sum_{i=1}^sx_i\tau(g_i)$ with $x_i\in\qwga$ (see [RW3, Lemma
 5]) and $\chi=(\infl)(\ti\chi)$, we compute
 $(\ti\tr\circ\defl)(x)(\ti\chi)=\ti\tr(\sum_i\ti x_i\tau(\ti
 g_i))(\ti\chi)=\sum_i\ol{\ti x_i}\,\ol{\ti g_i}\ti\chi(\ti
 g_i)=\sum_i\ol x_i\ol g_i\chi(g_i)$ and
 $(\defl\circ\tr)(x)(\ti\chi)=\defl(\tr\, x)(\ti\chi)=(\tr\,
 x)(\infl\ti\chi)=(\tr\, x)(\chi)=\sum_i\ol x_i\ol g_i\chi(g_i)$\,.
 Thus $(\ti\tr\circ\defl)(x)$ and $(\defl\circ\tr)(x)$ agree on
 irreducible characters $\ti\chi$ of $\ti G/\ti\Ga$. Since every
 irreducible character of $\ti G$ is a twist $\ti\chi\ot\ti\rho$
 with $\ti\rho$ of type W (compare [RW4, p.164, proof of 2.~of
 Lemma 4]), checking compatibility with W-twisting then gives
 $(\ti\tr\circ\defl)(x)=(\defl\circ\tr)(x)$\,.
 \sn The lemma is established.
 \Proposition{2.2} If $\tki\in T(\lwgi)$\,, then
 $\tki$ is in $\Ll(K_1(\lwgi))$. Moreover, there exists a $y\in(\lwgi)\mal$ so that $\Ll(y)=\tki$
 and $\deflab(y)=\la_{\kab/k}$\,, where $\kab$ is the fixed field of\, $[\gi,\gi]\,(\le
 H)$ and $G^\mr{ab}=\gi/[\gi,\gi]$.\Stop
 The proof depends on the commutative diagram
 $$\barr{cccccc}
 1+\fa_\wedge&\into&(\lwgi)\mal&\sr{\defl}{\onto}&(\lwgab)\mal&\\
 \sda&&\daz{\Ll}&&\daz{\Ll^\mr{ab}}&\\
 \tau(\fa_\wedge)&\into&T(\lwgi)&\sr{\defl}{\onto}&\lwgab&.\earr$$
 To understand the diagram, we recall that the canonical map
 $(\lwgi)\mal\to K_1(\lwgi)$ is surjective and
 that $(\lwgab)\mal=K_1(\lwgab)$ (see [CR, 40.31 and 45.12]). Whence we may regard $\Ll$
 as defined on $(\lwgi)\mal$.
 \sn In the diagram, $1+\fa_\wedge$ is the kernel of the upper defl.
 This map is surjective since $\fa_\wedge\subset\mr{rad}(\lwgi)$.
 The surjectivity of the lower defl is obvious. Finally, by [RW3,
 2b.~of Proposition 11] the left vertical map is surjective.
 The snake lemma shows that $\coker(\Ll)$ and $\coker(\Ll^\mr{ab})$ are
 isomorphic and that $\ker(\Ll)$ maps onto $\ker(\Ll^\mr{ab})$.
 Since $\tkab=\Ll^\mr{ab}(\la_{\kab/k})$ and $\defl(\tki)=\tkab$, the
 first assertion follows.
 \sn So $\Ll(y)=\tki$ for some $y\in(\lwgi)\mal$. Then
 $\defl(y)\me\la_{\kab/k}\in\ker(\Ll^\mr{ab})$, by the diagram, hence it
 equals $\defl(z)$ with $z\in\ker(\Ll)$. Thus $\defl(yz)=\la_{\kab/k}$
 and $\Ll(yz)=\tki$\,. \sn The proof of Proposition 2.2 is complete.
 \bn To continue the discussion of the integrality of $\tki$,
 we again pick a central open subgroup $\Ga$ in $\gi$
 and representatives $g_1,\ldots, g_s$ in $\gi$ of the conjugacy
 classes of the group $\gi/\Ga$; so the images $\tau(g_i)$ of
 $g_i$ in $T(\lwgi)$ constitute a $\lwga$-basis of $T(\lwgi)$
 (see [RW3, Lemma 5$_\wedge$]).
 We will employ the following notation\,:
 \ben\item[]{\em given $t\in T(\qwgi)$, there is a unique function
 $\tga:\gi\to\qwga$ such that \bit\item[]$\tga(g)=\tga(g')$ if $g$
 and $g'$ are conjugate in $\gi$\,,
 \item[]$\tga(\ga g)=\ga\tga(g)$ for $\ga\in\Ga$\,,
 \item[]$t=\sum_{i=1}^r\tga(g_i\me)\tau(g_i)$ with the $g_i$ as
 above.\eit}\een
 This is easily checked. We call $\tga(g\me)$ the {\em
 coefficient}
 of $t$ at $g$.
 In particular, [\,$\tki\in T(\lwgi)\iff
 (\tki|\Ga)(g\me)\in\lwga$ for all $g\in\gi$\,]\,.
 \mn We now compute $(\tki|\Ga)(g\me)$. From $\tr(\tki)=\LL(\lki)$ we
 obtain $$\sum_{i=1}^r\ol{(\tki|\Ga)(g_i\me)}\,\ol
 g_i\chi_j(g_i)=\fra1l\log\fra{\lki(\chi_j)^l}{\Psi\lki(\plcj)}$$
 for the irreducible characters $\chi_j,\,1\le j\le s$\,, of
 $\gi/\Ga$. This equation enables us to perform a Fourier
 inversion in order to isolate the $\ol{(\tki|\Ga)(g_i\me)}$. To that end,
 let $h_i$ denote the order of the conjugacy class of $(g_i\mod\Ga)$
 and denote by $M\,,\,M_1$ the $s\times s$-matrices
 $(\chi_i(g_j))_{i,j}\,,\,(h_i\chi_j(g_i\me))_{i,j}$, respectively. The
 orthogonality relations yield $M\cdot M_1=[\gi:\Ga]\cdot{\ul 1}$
 whence $M_1\cdot M=[\gi:\Ga]\cdot\ul{1}$ as well, and we arrive first
 at
 $$M\cdot\left(\barr{c}\vdots\\ \ol{(\tki|\Ga)(g_i\me)}\,\,\ol g_i\\
 \vdots\earr\right)=\left(\barr{c}\vdots\\
 \fra1l\log\fra{\lki(\chi_j)^l}{\Psi\lki(\plcj)}\\
 \vdots\earr\right)$$
 and then at $$[\gi:\Ga]\ol{(\tki|\Ga)(g_i\me)}\,\ol
 g_i=\sum_{j=1}^sh_i\chi_j(g_i\me)\fra1l\log\fra{\lki(\chi_j)^l}{\Psi\lki(\plcj)}\,,\quad(1\le
 i\le s)\,.$$
 So we have
 \Proposition{2.3}
 $\ol{(\tki|\Ga)(g_i\me)}=\fra1{[\gi:\Ga]}\sum_{j=1}^sh_i\chi_j(g_i\me)\fra1l\log\fra{\lki(\chi_j)^l}
 {\Psi\lki(\plcj)}\ol g_i\me$
 \Stop
 Observe that the right hand side is in $\qwgak$ (by Galois invariance); it is in $\lwgak$
 for all $g_i$ precisely when $\tki\in T(\lwgi)$.
 Recall also that, up to W-twists, the $\chi_j$ are a full set
 of representatives of all irreducible characters of $\gi$.
 \Proposition{2.4} \ben\item If $\tki\in T(\lwgi)$, then
 there exists a torsion element  $$w\in\HOM(\rlgi,(\lcwgak)\mal) \
 so\
 that \ \deflab(w)=1 \ and \ w\cdot\lki\in\Det\, K_1(\lwgi)\,.$$
 \item If $w\in\HOM(\rlgi,(\lcwgak)\mal)$ is torsion and $w\cdot\lki\in\Det\,K_1(\lwgi)$\,,
 then $\tki\in
 T(\lwgi)$\,.
 \item There is at most one torsion element  $w\in\HOM(\rlgi,(\lcwgak)\mal)$  so
 that $\deflab(w)$ $=1$ and $w\cdot\lki\in\Det\, K_1(\lwgi)$\,.\een\Stop
 For 1., take a $y$ with $\Ll(y)=\tki$ from Proposition 2.2 and define $w$
 by $w\cdot\lki=\Det(y)$. Applying $\defl=\deflab$, we get
 $$\defl(w)L_{\kab/k}=\Det(\defl\,y)=\Det(\la_{\kab/k})=L_{\kab/k}\,,$$
 hence $\defl(w)=1$. Next apply $\LL$ of {\sc Fact 2} and get
 $$\LL(w)+\tr(\tki)=\LL(\Det\,y)=\tr(\Ll\,y)=\tr(\tki)\,,$$ hence
 $\LL(w)=0$. By the argument in [RW3, \S6, just before the
 Corollary of Theorem B$_\wedge$], showing that $\Det\,z$ is
 torsion, it follows that $w$ is torsion.
 \sn 2.~follows from the ${\sss{\wedge}}$-version of {\sc Fact 2} and
 $\LL(w)=0$.
 \sn For 3., suppose $w_1,w_2$ are both as above, hence
 $w_i\lki=\Det(y_i)$ with $y_i\in(\lwgi)\mal$ and $\defl(w_i)=1$.
 Then $w_1\me w_2=\Det(y_1\me y_2)$ and $\defl(w_1\me w_2)=1$.
 Therefore, $\Det(\defl(y_1\me y_2))=1$, since $\Det$ commutes
 with $\defl$, from which $\defl(y_1\me y_2)=1$ follows because
 Det is an isomorphism by $SK_1(\La_\wedge G^\mr{ab})=1$. So
 $y_1\me y_2\in 1+\fa_\wedge$, by the diagram in the proof of
 Proposition 2.2, and thus $w_1\me w_2\in\Det(1+\fa_\wedge)$ is
 torsion and therefore 1 by [RW3, Lemma 12].
 \sn The proof is finished.
 \remark  It may be worth mentioning that in Fröhlich's work on
 Galois module structure logarithmic and torsion congruences
 appear as well: see [Fr, IV,\S\S4,5,6]. Here, the logarithmic
 congruences are taken care of by the Davenport-Hasse formula
 [Fr, p.179].
 \Section{3}{The torsion congruence}
 In this section $S$ is sufficiently large. We continue to assume $\mu=0$ for the pro-$l$ extension $\Ki/k$ and,
 under the hypothesis $\tki\in T(\lwgi)$, exhibit the additional
 congruences which are equivalent to
 $\lki\in \Det\,K_1(\lwgi)$ or, by {\sc Fact 3},
 to $\lki\in\Det\,K_1(\lbgi)$.
 \mn Assume first that there is an abelian subgroup $G'$ of index $l$ in
 the non-abelian pro-$l$ group $\gi$. Denote by $k'$ the fixed
 field of $G'$ and by $\ver$ the transfer map $\giab\to G'$.
 Setting $A=\gi/G'$, consider the conjugation action of
 $A$ on $G'$ and $\La_\wedge\gs$, and let ${\cal{T}'}$ be the image of the $A$-trace map
 on
 $\La_\wedge\gs$. Thus ${\cal{T}'}$ is an ideal in the ring $(\La_\wedge\gs)^A$ of
 $A$-fixed points of $\La_\wedge\gs$.
 \Lemma{3.1} $\la_{\Ki/k'}$ is $A$-invariant.\Stop
 This is because
 $\Det:K_1(\lbgs)\to\HOM(R_l(G'),\La_\bullet\clo(\Ga_{k'})\mal)$
 is an $A$-equivariant monomorphism (see [RW4, 2\,nd paragraph on
 p.159]) and
 $\Det\,\la_{\Ki/k'}=\lks\,,\,L_{\Ki/k'}^a(\chi')=\lks({\chi'}^{a\me})$
 $=\lki(\ind_{G'}^\gi({\chi'}^{a\me}))=\lki(\ind_{G'}^\gi(\chi'))=
 \lks(\chi')$ for all $a\in A$ and all $\chi'\in R_l(G')$.
 \Proposition{3.2} In the above situation and with $w$ as in Proposition 2.4, the following are equivalent
 \ben\item $w=1$, i.e., $\lki\in \Det\,K_1(\lwgi)$\,,
 \item $\ver(\la_{\kab/k})\equiv\la_{\Ki/k'}\mod {\cal{T}'}$\,,
 \item $\fra1{[G':\Ga]}\sum_{\chi'}\fra{\lks(\chi')}{\Psi\lki(\cver)}\equiv 1\mod l\cdot
 \lcwgak$\,, with the sum ranging over the $\ql\clo$-irreducible characters
 $\chi'$ of $G'/\Ga$ where $\Ga$ is a central open subgroup of $\gi$.\een\Stop
 \proof Note first that $w=1$ if, and only if,
 $\res_\gi^{G'}w=1$\,: If $\chi$ is an irreducible character of
 $\gi$, then either $\chi=\infl_{G^\mr{ab}}^\gi\al$ is
 inflated from an abelian character $\al$ of $G^\mr{ab}$ or
 $\chi=\ind_{G'}^\gi(\chi')$ is induced from an abelian character
 $\chi'$ of $G'$. Therefore, either
 $w(\chi)=w(\infl\,\al)=(\defl\,w)(\al)=1$ or
 $w(\chi)=w(\ind\,\chi')=(\res\,w)(\chi')=1$\,.
 \mn Write $w\cdot\lki=\Det\,y$ with $y\in(\lwgi)\mal$ as in the proof of Proposition 2.4, so that
 $\deflab\,y=\la_{\kab/k}$. Since $\res_\gi^{G'}$ commutes with
 Det and $\res_\gi^\gs\lki=\lks$ (see [RW2, Lemma 9 and p.563]),
 we have
 $$\res_\gi^{G'}w=\Det(\fra{\res_\gi^{G'}y}{\la_{\Ki/k'}})\,.$$
 Since Det is injective on $(\La_\wedge G')\mal$, it follows that
 $\fra{\res_\gi^{G'}y}{\la_{\Ki/k'}}$ is a torsion element.
 \mn By means of the commutative square
 \par\bmp{6cm}$\barr{cccc}(\lwgi)\mal&\onto&K_1(\lwgi)&\\
 \daz{N}&&\daz{\res_\gi^\gs}&\\
 (\La_\wedge\gs)\mal&=&K_1(\La_\wedge\gs)&\earr$
 \emp\hsp{10}\bmp{7cm} we write $\res_\gi^{G'}y=N(y)$, the
 determinant of right multiplication of $y$ on the left
 $\La_\wedge\gs$-module $\lwgi$ (compare [RW3, proof of Lemma 12]).\emp
 \mn Denoting the composite map
 $\gi\sr{\defl}{\lto}G^\mr{ab}\sr{\ver}{\lto}G'$ also by
 ver, we obtain from $\defl\,y=\la_{\kab/k}$ the trivial equation
 $$\fra{\res_\gi^{G'}y}{\la_{\Ki/k'}}=\fra{N(y)}{\ver(y)}\cdot\fra{\ver(\la_{\kab/k})}{\la_{\Ki/k'}}\
 .$$ By a congruence due to C.T.C.~Wall, \ $N(y)\equiv\ver(y)\mod
 {\cal{T}'}$ \ with ${\cal{T}'}$ the $A$-trace ideal
 of the $A$-action on $\La_\wedge G'$ (see [RW3, proof of Lemma 12]).
 Since both, $N(y)$ and $\ver(y)$, are units in $\La_\wedge G'$
 fixed by $A$, it follows that
 $$\fra{\res_\gi^{G'}y}{\la_{\Ki/k'}}\equiv
 \fra{\ver(\la_{\kab/k})}{\la_{\Ki/k'}}\mod {\cal{T}'}\,.$$
 The equivalence of 1.~and 2.~now follows from the last two
 paragraphs if we can show for torsion units $e$ of $\La_\wedge
 G'$ that $$e\equiv1\mod {\cal{T}'}\iff e=1\,.$$ But this
 has already been shown in the fifth paragraph of the proof of
 [RW3, Lemma 12]\,: the argument
 $$N(x)\ \mr{torsion}\ \ \&\ \ N(x)\equiv1\mod\cT\implies N(x)=1$$ works with $N(x)$ replaced by $e$.
 \mn We next turn to the equivalence of 2.~and 3.~and first prove,
 for irreducible characters $\chi'$ of $G'/\Ga$,
 $$\Det(\fra{\la_{\Ki/k'}}{\ver(\la_{\kab/k})})(\chi')=\fra{\lks(\chi')}{\Psi(\lki(\chi'\circ\ver))}\
 .$$
 With $\la_{\kab/k}=\defl\,y$ as before, hence
 $\ver(\la_{\kab/k})=\ver(y)$,
 write $y=\sum_{x}y_xx$ with $y_x\in\lwga$ and
 $\{x\}$ a set of representatives of the elements of $\gi/\Ga$ in
 $\gi$.
 Since ver induces $\Psi$ on $\Ga$, we obtain
 $\ver(y)=\sum_x\Psi(y_x)\ver(x)$\,, so
 $$(\Det(\ver\,y))(\chi')=\sum_x\ol{\Psi(y_x)}\ol
 x^l\chi'(\ver\,x)$$ because $\ol{\ver(x)}=\ol x^l$ (compare the first
 displayed formula after {\sc Fact 1}).
 \sn On the other hand, $(\Det\, y)(\cver)=\sum_x\ol y_x\ol
 x(\cver)(x)$\,, whence  $\Psi(\Det\, y)(\cver)=\sum_x\Psi(\ol y_x)\ol x^l(\cver)(x)$\,.
 Thus $(\Det(\ver\,y))(\chi')=\Psi((\Det\,y) (\chi'\circ\ver))$ and we
 conclude {\small
 $$(\Det\fra{\la_{\Ki/k'}}{\ver(y)})(\chi')=\fra{\lks(\chi')}{\Psi((\Det\,
 y)(\cver))}=\fra{\lks(\chi')}{\Psi((w\cdot\lki)(\cver))}=\fra{\lks(\chi')}{\Psi((\lki)(\cver))}$$ }
 as $\deflab\,w=1$.
 \sn Now write
 $\fra{\la_{\Ki/k'}}{\ver(\la_{\kab/k})}=\sum_{j=1}^sa_jg_j$ with
 $a_j\in\La_\wedge\Ga$ and $1=g_1,\ldots,g_s$ a set of representatives of the elements of $\gs/\Ga$ in
 $\gs$, of which the first $r$ are precisely those in the centre $Z(\gi)$ of $\gi$.
 Then
 $$\Det(\fra{\la_{\Ki/k'}}{\ver(\la_{\kab/k})})(\chi')=\sum_j\ol
 a_j\ol g_j\chi'(g_j)\ ,$$ and by Fourier inversion, in view of
 the above paragraph,
 $$\ol{
 a_jg_j}=\fra1{[G':\Ga]}\sum_{\chi'}\chi'(g_j\me)\fra{\lks(\chi')}{\Psi(\lki(\chi'\circ\ver))}\
 .$$
 Since, by Lemma 3.1, $\fra{\la_{\Ki/k'}}{\ver(\la_{\kab/k})}$ is
 $A$-invariant, the coefficients $a_j$ are constant on orbits of
 the $A$-action on $G'$, hence
 $\fra{\la_{\Ki/k'}}{\ver(\la_{\kab/k})}\equiv\sum_{j=1}^ra_jg_j\mod\cT$.
 Note that by the fourth paragraph \footnote{with care taken to choose the coset
 representatives $\{b\}$ closed under conjugation by $a$; this is possible since
 $(b\Ga)^a=b\Ga$
 implies $b^a=b$, as $b^{a-1}\in\Ga\cap[\gi,\gi]=1$} of the proof of [RW3, Lemma 12]
 $g_1+{\cal{T}'},\ldots,g_r+{\cal{T}'}$
 is a $\La_\wedge\Ga/l\La_\wedge\Ga$-basis of $(\La_\wedge
 G')^A/{\cal{T}'}$\,.
 \sn The fifth paragraph of that proof
 shows that $\fra{\la_{\Ki/k'}}{\ver(\la_{\kab/k})}\equiv\ze z\mod
 {\cal{T}'}$\,, with a root $\ze$ of unity and a $z\in Z(\gi)$ of
 finite order. If $\ze
 z\equiv\sum_{j=1}^ra_jg_j\mod{\cal{T}'}$\,, then there is
 a unique $j_0$ so that $\ze z\equiv a_{j_0}g_{j_0}$ and all
 $a_jg_j\equiv 0$ for $j\neq j_0$ (this follows by writing $z=\ga
 g_{j_0}$ with $\ga\in\Ga$). Hence, 2.~is equivalent to
 $a_1\equiv1\mod l\La_\wedge\Ga$ (since this implies $j_0=1)$. But
 the latter is equivalent to $\ol a_1\equiv1\mod l\La_\wedge\ol\Ga$\,,
 thus to $\ol a_1\equiv1\mod l\La_\wedge\Ga_k$ as
 $\La_\wedge\ol\Ga\cap l\La_\wedge\Ga_k=l\La_\wedge\ol\Ga$ where
 $\La_\wedge\Ga_k$ is a free $\La_\wedge\ol\Ga$-module.
 \remark The implication $[\,2.\implies
 3.\,]$ does not need any hypothesis as it depends on the Fourier
 inversion step in the above proof. Thus 3.~can be
 restated as
 $$\fra1{[G':\Ga]}\sum_{\chi'}\chi'(z\me)\Big(\fra{L_{\Ki/k'}(\chi')}{\Psi\lki(\chi'\circ\ver)}-1\Big)\equiv0\mod
 l\cdot\lcwgak\quad (\forall\,z\in Z(\gi))\,,$$
 since the coefficients of
 $\fra{\la_{\Ki/k'}}{\ver(\la_{\kab/k})}-1$ at the central
 elements are divisible by $l$.
 \mn We finally free ourselves from the assumption that $\gi$ has
 an abelian subgroup of index $l$.
 \theorem Assume $\mu=0$ for $\Ki/k$ and $\tki\in T(\lwgi)$. Then
 $\lki\in\Det\, K_1(\lwgi)$ provided that $L_{F/f}\in\Det\,K_1(\La_\wedge G_{F/f})$ for all intermediate Galois
 extensions $F/f$ in $\Ki/k$ such that
 \ben\item $G_{F/f}$ has an abelian subgroup of
 index $l$\,,
 \item $[\Ki:F]$ is finite and $f$ is fixed by the centre
 $Z(\gi)$ of $\gi$\,.\een\Stop
 The proof is by induction on $[\gi:Z(\gi)]$.\sn By 1.~of Proposition 2.4,
 there exists a torsion $w\in\HOM(\rlgi,(\lcwgak)\mal)$ so that $\deflab
 w=1$ and $w\cdot\lki\in\Det\,K_1(\lwgi)$\,. It suffices to
 show that $w=1$.
 \sn Let $G'$ be any subgroup of index $l$ in $\gi$ containing
 $Z(\gi)$ and let $ k'$ be its fixed field. Note that all
 intermediate extensions of $\Ki/k'$ satisfying conditions 1.~and
 2.~relative to $\Ki/k'$ also satisfy 1.~and 2.~relative to $\Ki/k$ because
 $Z(G')\supset Z(\gi)$. So the induction hypothesis implies
 $L_{\Ki/k'}\in \Det\,K_1(\La_\wedge G')$.
 \sn
 Consider the commutative diagram
 $$\barr{ccccc}{\cal{H}}(\Ki/k)&\sr{d}{\lto}&{\cal{H}}(\ti K/k)&\sr{q}{\lto}&{\cal{H}}
 (\kab/k)\\
 \daz{\res_\gi^\gs}&&\daz{\res_{\ti G}^{{G'}^\mr{ab}}}&&\\
 {\cal{H}}(\Ki/k')&\sr{d'}{\lto}&{\cal{H}}(\ti K/k')&&\earr$$
 with $\ti G=G_{\ti K/k}=\gi/[G',G']$ and $\ti K$ the fixed field of $[G',G']$. Here, ${\cal{H}}(F/f)$ abbreviates
 $\Hom(R_l(G_{F/f}),(\La_\wedge\clo\Ga_f)\mal)$ and all horizontal maps are deflations
 \footnote{this is the analogue of the commutative diagram in [RW3, proof of
 Lemma 12]}.
 \sn Now,
 $$(dw)L_{\ti K/k}=d(w\lki)\in\Det\,K_1(\La_\wedge\ti G)\,,\
 q(dw)=\deflab\,w=1$$
 (by $[\gs,\gs]\subset[\gi,\gi]$). Since $\ti G$ has the
 abelian subgroup ${G'}^\mr{ab}$ of index $l$, $dw=1$ follows by hypothesis 1.~and 2.
 \sn Thus $d'(\res_\gi^{G'}w)=\res_{\ti G}^{{G'}^\mr{ab}}(dw)=\res_{\ti
 G}^{{G'}^\mr{ab}}1=1$ and
 $(\res_\gi^{G'}w)L_{\Ki/k'}\in\Det\,K_1(\La_\wedge G')$
 (see [RW2, 2.~of Proposition 12]). Hence, by the above,
 $$\res_\gi^{G'}w\in\Det\,K_1(\La_\wedge G')\subset\HOM(R_l(G'),(\La\clo_\wedge\Ga_{k'})\mal)\,,$$ by
 [RW3, 1.~of Proposition 11]
 We can now invoke 3.~of Proposition 2.4 to conclude that $\res_\gi^{G'}w=1$.
  \sn We now show $w(\chi)=1$ for every irreducible $\chi$.  If
 $\chi$ is abelian, then $\deflab w=1$ takes care of this. If
 $\chi$ is non-abelian, then it is induced from a character
 $\chi'$ of a subgroup $G'$ of index $l$ in $\gi$ which contains
 $Z(\gi)$ (see [CR, 11.2]). Hence,
 $w(\chi)=(\res_\gi^{G'}w)(\chi')=1\,,$
 establishing the theorem.
 \mn{\sc Remarks.}\ben\item Combining the theorem with {\sc Fact 3} we even get $\lki\in
 \Det\,K_1(\lbgi)$, hence $\Ki/k$ has a non-abelian pseudomeasure.
 Moreover, in all three statements in Proposition 3.2
 $\pht{.}_\wedge$ can be replaced by\,$\pht{.}_\bullet$\,; for
 3.~this is because
 $\lbgak\cap l\lwgak=l\lbgak$\,; for 2.~note that
 ${\cal{T}'}\cap\lbgs$ is the image of the $A$-trace map on
 $\lbgs$.
 \item The torsion congruence in Proposition 3.2 is
 somehow reminiscent of special value congruences which have been
 derived in [Ty1].\een
 \Section{4}{Minimal non-abelian $\gi$}
 In this section $\gi$ is a non-abelian pro-$l$
 group which has an abelian subgroup $\gs$ of index $l$. Further, $A$ is short for the cyclic group
 $\gi/\gs$ of order $l$ and $a\mod\gs$ is a generator. Finally, we fix an irreducible character
 $\om\in\rlgi$ which is trivial on $\gs$
 but not on $a$ and a central open subgroup $\Ga$
 of $\gi$ (so $\Ga\le Z(\gi)\le\gs\le\gi$). \sn We continue to
 assume $\mu=0$ for $\Ki/k$ and $S$ to be sufficiently large.
 \Lemma{4.1} \ben\item If $\chi'$ is an irreducible character of
 $\gs$, then
 $$\ind_\gs^\gi(\plcs)-\psi_l(\ind_\gs^\gi\chi')=\sum_{i=0}^{l-1}(\csver)\om^i-l(\csver)$$
 with {\rm ver} the transfer $\gi\to\giab\to\gs$\,,
 \item $\fra{[\gi:Z(\gi)]}{l\cdot|[\gi,\gi]|}\in\Z$\,.\een\Stop
 The first assertion results on evaluating both sides on an
 element $g\in\gi$ (recall $\ver(g)=\prod_{i=0}^{l-1}g^{a^i}$ or
 $=g^l$ according as $g\in\gs$ or not, and $\im(\ver)\subset Z(\gi)$).
 \sn For the second assertion we use induction on $|[\gi,\gi]|$. By
 $[\gi,\gi]=[\gi,\gs]$ we can pick
 a central commutator $z=[a,g']$ of order $l$; we set $\ti G=\gi/\group{z}$. If
 $|[\gi,\gi]|=l$, then $\group{z}=[\gi,\gi]$ and the fraction in question is
 $\fra{[\gi:Z(\gi)]}{l^2}$, which is integral as $\gi$ is
 non-abelian. If however $|[\gi,\gi]|>l$, then $\ti G$ is non-abelian, but $[\ti G,\ti
 G]=\widetilde{[\gi,\gi]}$ and $Z(\ti G)\varsupsetneq\widetilde{Z(\gi)}$\,; the
 latter by $\ti g'\in Z(\ti G)\setminus\widetilde{Z(\gi)}$\,.
 \sn Summing up and setting $l^\nu=[Z(\ti G):\widetilde{Z(\gi)}]$\,, we obtain
 $\Z\ni\frac{[\ti G:Z(\ti G)]}{l\cdot|[\ti G,\ti
 G]|}=\frac{[\gi:Z(\gi)]}{|[\gi,\gi]|\cdot l^\nu}$\,, which finishes the
 proof of 2.~and thus of the lemma.
 \mn We next turn to the computation of the coefficient
 $\tkiga(g\me)$ at $g\in\gi$ of the logarithmic pseudomeasure
 $\tki$. To ease notation, set $m_\chi=\LL(\lki)(\chi)=\mchi$
 whenever $\chi$ is a character of $\gi$ which is trivial on
 $\Ga$. Obviously, $m_{\chi_1+\chi_2}=m_{\chi_1}+m_{\chi_2}$ and
 $m_{\infl_{\ti G}^\gi\ti\chi}=\ti m_{\ti\chi}$ for characters
 $\ti\chi$ of $\ti G=\gi/N$, $N\lhd\gi$ finite, with $\ti\Ga=\Ga
 N/N\subset\ker\ti\chi$\,.
 \Proposition{4.2} $\tkiga(g\me)\in\lwga$ \ for all $g\in\gi\setminus
 Z(\gi)$\,.\Stop
 By Proposition 2.3, $\ol{\tkiga(g\me)}=\giga\sum_\chi
 h_g\chi(g\me)m_\chi\ol g\me$\,, with the sum ranging over the irreducible
 $\chi\in\rlgi$ which have $\Ga$ in their kernel. We need to show
 $$\giga\sum_\chi h_g\chi(g\me)m_\chi\in\lwgak\ \ \mr{for}\ \ g\in\gi\setminus
 Z(\gi)\,.$$
 \sn If $g\notin\gs$, then the centralizer of $g$ in $\gi$ is
 $Z(\gi)\cdot\group{g}$, hence $h_g=\fra1l[\gi:Z(\gi)]$. Moreover,
 $\chi(g)=0$ for every non-abelian $\chi$, since these are induced
 from $\gs$. Therefore, with $N=[\gi,\gi]$,
 $$\barr{l}\giga\sum_\chi
 h_g\chi(g\me)m_\chi=\fra1{l\cdot[Z(\gi):\Ga]}\sum_{\chi(1)=1}\chi(g\me)m_\chi=\fra{[\gi:N\cdot\Ga]}
 {l\cdot[Z(\gi):\Ga]}\ol{\tkaga(\ti g\me)}\ol g\,.\earr$$
 This is in $\lwgak$ since $N\cap\Ga=1$ and $\fra{[\gi:N\Ga]}{l\cdot[Z(\gi):\Ga]}=
 \fra{[\gi:\Ga]}{l\cdot [Z(\gi):\Ga]|N|}=\fra{[\gi:Z(\gi)]}{l\cdot|N|}\in\Z\,,$ by 2.\linebreak of Lemma 4.1 and
 $t_{\kab/k}=\Ll(\la_{\kab/k})\in T(\La_\wedge G^\mr{ab})$\,.
 \sn Next, let $g\in\gs$. Then
 $$\giga\sum_\chi
 h_g\chi(g\me)m_\chi=\fra{h_g}{[\gi:\Ga]}\Big(\sum_{\chi(1)=1}\chi(g\me)m_\chi+\sum_{\chi(1)=l}\chi(g\me)m_\chi\,\Big)\,.$$
 If $\chi(1)=1$ set $\chi'=\res_\gi^\gs\chi$. Then $\chi'$ is
 trivial on $[\gi,\gi]$. We denote the sum over all such $\chi'$
 by $\sum_1$ and obtain
 $$\barr{l}\sum_{\chi(1)=1}\chi(g\me)m_\chi=\sum_1\sum_{i=0}^{l-1}(\chi\om^i)(g\me)m_{\chi\om^i}\\[1mm]
 =\sum_1\chi'(g\me)\fra1l\log\prod_{i=0}^{l-1}\fra{\lki(\chi\om^i)^l}{\Psi\lki(\psi_l(\chi\om^i))}=
 \sum_1\chi'(g\me)\fra1l\log\fra{\lks(\chi')^l}{\Psi\lki(\psi_l(\ind_\gs^\gi\chi'))}\,.\earr$$
 Above, we have used $\om(g)=1$ and
 $\lki(\ind_\gs^\gi\chi')=\lks(\chi')$\,.
 \sn If $\chi(1)=l$, then $\chi=\ind_\gs^\gi\chi'$ with an abelian
 character $\chi'$ of $\gs/\Ga$ which is non-trivial on
 $[\gi,\gi]$. We denote the sum over all such $\chi'$ by $\sum_2$
 and obtain
 $$\barr{l}\sum_{\chi(1)=l}\chi(g\me)m_\chi=\fra1l\sum_2(\ind_\gs^\gi\chi')(g\me)\fra1l\log\fra{\lks(\chi')^l}
 {\Psi\lki(\psi_l(\ind_\gs^\gi\chi'))}\\[1mm]
 =\fra1l\sum_2\sum_{i=0}^{l-1}\chi'(g^{-a^i})\fra1l\log\fra{\lks({\chi'}^{a^i})^l}{\Psi\lki(\psi_l
 (\ind_\gs^\gi{\chi'}^{a^i}))} =
 \sum_2\chi'(g\me)\fra1l\log\fra{\lks(\chi')^l}{\Psi\lki(\psi_l(\ind_\gs^\gi\chi'))}\,,\earr$$
 where we have used that $\chi'$ and ${\chi'}^{a^i}$ induce the
 same $\chi$.
 \sn Collecting everything so far, we see that
 \nf{$\star$}
 {\giga\sum_\chi h_g\chi(g\me)m_\chi=\fra{h_g}{[\gi:\Ga]}\sum_{\chi'}\chi'(g\me)\fra1l
 \log\fra{\lks(\chi')^l}{\Psi\lki(\psi_l(\ind_\gs^\gi\chi'))}\, .}
 Multiplying the term in log by \
 $1=\fra{\Psi\lki(\ind_\gs^\gi\plcs)}{\Psi\lks(\plcs)}$
 \ gives
 \nf{$\sharp$}{\hsp{-3}\barr{l}\fra{h_g}{l}\ol{(t_{\Ki/k'}|\Ga)}(g\me)\ol
 g+\fra{h_g}{[\gi:\Ga]}\sum_{\chi'}\chi'(g\me)\fra1l\log\Psi\lki(\indpsi)\earr}
 by 1.~of Lemma 4.1. Since $g\not\in Z(\gi)$, $h_g=l$ and
 $\ol{(t_{\Ki/k'}|\Ga)(g\me)}\ol g$
 is integral by $t_{\Ki/k'}=\Ll(\la_{\Ki/k'})\in T(\La_\wedge G')$. With respect to the second summand we observe that
 multiplying $\chi'$ by an abelian character $\theta$ of
 $\gs/\ver(G^\mr{ab})$ does not change the character
 $\sum_{i=0}^{l-1}(\chi'\circ\ver)\om^i-$ $l(\chi'\circ\ver)$. This brings
 the sum $\sum_\theta\theta(g\me)$ into the second summand and makes it
 vanish, if $g\notin \ver(G^\mr{ab})$.
 \sn The proof of the proposition is complete since $\ver(G^\mr{ab})\subset Z(\gi)$.
 \corollary
 $\ol{(\tki|\Ga)(g\me)}=\fra{h_g}{l}\ol{(t_{\Ki/k'}|\Ga)(g\me)}$ \
 for all $g\in G'$ which are not in the image of $\ver$.\Stop
 \hsp{-3.6} To continue with computing the coefficients
 $(\tki|\Ga)(z\me)$ for $z\in Z(\gi)$ we add the assumption
 $[\gi,\gi]\subset Z(\gi)$.
 \Lemma{4.3} If $[\gi,\gi]\subset Z(\gi)$, then, with $a\mod G'$ generating
 the cyclic group $\gi/G'$ of order $l$,
 \ben\item $\forall\,g'\in G'\ \exists\,z\in[\gi,\gi]: \ (g')^{a^i}=g'z^i\,,\,z^l=1$
 \item
 ${\chi'}^{\hat A}\df\prod_{i=0}^{l-1}{\chi'}^{a^i}={\chi'}^l$ for
 all irreducible characters $\chi'$ of $G'/\Ga$
 \item $x\in\cT\implies x^{l^v}\in l^v\cT$\ .
 \een\Stop
 Indeed, $z\df (g')^{a-1}\in[\gi,\gi]$, so $z^l=(g')^{a^l-1}=1$\,,
 and $(g')^{\hat A}={g'}^l$ implies ${\chi'}^{\hat A}={\chi'}^l$.
 \sn For 3., take a set $\{b\}$ of representatives in $G'$ of the orbits
 of $G'/\Ga$ under conjugation by $A$. Then $\{\mr{tr}_A(b)\}$ is
 a $\La_\wedge\Ga$-basis of $\cT$, so we can write
 $x=\sum_bx_b\mr{tr}_A(b)$. We show that $x^l\in l\cT$\,: \sn Since
 $x^l\equiv\sum_bx_b^l\mr{tr}_A(b)^l\mod l\cT$, this follows from
 $\mr{tr}_A(b)^l\in l\cT$. To see the latter, write
 $b^a=bz\,,\,\hat z=1+z+\ldots+z^{l-1}$\,, so
 \ $\mr{tr}_A(b)^l=(b\hat z)^l=b^l{\hat z}^l=l^{l-1}b^l\hat
 z=l^{l-2}(lb^l\hat z)\in l^{l-2}\cT$ \ (by $(b^l)^a=b^lz^l=b^l$), as
 required since $l$ is odd.
 \sn 3.~now follows by induction on $v$.
 \corollary If $x\in \cT$, then for every $z\in Z(\gi)$,
 $$\fra{1}{[G':\Ga]}\sum_{\chi'}\chi'(z\me)\fra1l\log(\,\Det(1+x)(\chi')\,)\in\lcwgak\,.$$\Stop
 In fact, as in the remark following Proposition 3.2, we have (\,with
 $x=\fra{\la_{\Ki/k'}}{\ver(\la_{\kab/k})}-1\,)$
 \nf{$\Diamond$}{\fra{1}{[G':\Ga]}\sum_{\chi'}\chi'(z\me)(\Det\,x)(\chi')\equiv0\mod
 l\lcwgak\,.} Now, since
 $\Det(\sum_bx_b\mr{tr}_A(b))(\chi')=\sum_{b,i}\ol x_b\ol
 b\chi'(b^{a^i})=l\sum_b\ol x_b\ol b+\sum_{b,i}\ol x_b\ol
 b(\chi'(b^{a^i})-1)$ is divisible by $\ze-1$ for some $l$-power
 root $\ze$ of unity, we see that the logarithmic series
 $\log(1+\Det\,x)(\chi')$ converges in $\qwgak$, hence
 $$\log(\Det(1+x)(\chi'))=\log(\,1+(\Det\,x)(\chi')\,)=\sum_{\nu=1}^\infty(-1)^{\nu-1}\fra1\nu(\Det\,x^\nu)(\chi')\,.$$
 So we obtain
 $$\barr{l}\fra{1}{[G':\Ga]}\sum_{\chi'}\chi'(z\me)\fra1l\log(\,\Det(1+x)(\chi')\,)=
 \sum_{\nu=1}^\infty(-1)^{\nu-1}\fra1{l\nu}
 \fra1{[G':\Ga]}\sum_{\chi'}\chi'(z\me)(\Det\,x^\nu)(\chi')\\
 =\ol z\sum_{\nu=1}^\infty(-1)^{\nu-1}\fra{1}{l\nu}\ol{(x^\nu|\Ga)(z\me)}\earr$$
 by Proposition 2.3 (on identifying $\La_\wedge G'$ and
 $T(\La_\wedge G'))$.
 \sn Writing $\nu=rl^v$ with $l\nmid r$, then the above $\nu$-th
 summand is $\pm\fra1{rl^{v+1}}\ol{((x^r)^{l^v}|\Ga)(z\me)}$ with
 $x^r\in\cT$, hence $(x^r)^{l^v}\in l^v\cT$ by 3.~of Lemma 4.3.
 Thus the corollary follows from ($\Diamond)$.
 \Proposition{4.4} Assume\ben\item $\gi$ has an abelian subgroup
 $G'$ of index $l$\,,
 \item $[\gi,\gi]\subset Z(\gi)$\,,
 \item
 $\fra{\la_{\Ki/k'}}{\ver(\la_{\kab/k})}\equiv1\mod\cT$\,.\een
 Then $\tki\in T(\lwgi)$\ \footnote{and thus $\lki\in\Det\,K_1(\lbgi)$, by Proposition 3.2} .\Stop
 Because of Proposition 4.2 the proof only requires checking
 integrality of the coefficients $(\tki|\Ga)(z\me)$ for $z\in
 Z(\gi)$. Invoking $(\star)$ and $(\sharp)$ of that proposition we see that
 $$\barr{l}\ol{(\tki|\Ga)(z\me)}\ol z =\fra1l\ol{(t_{\Ki/k'}|\Ga)(z\me)}\ol
 z\ +\\ \quad+\giga\sum_{\chi'}\chi'(z\me)
 \fra1l\log\Big(
 \Psi\lki(\indpsi)\Big)\,.\earr$$
 By  $\ind_\gs^\gi(\res^\gs_\gi(\cver))=\sum_{i=0}^{l-1}(\cver)\om^i$
 and $\res_\gi^\gs(\cver)=\chisa$\,, the second summand equals
 $\giga\sum_{\chi'}(z\me)\fra1l\log\Psi\fra{\lks(\chisa)}{\lki(\chi'\circ\ver)^l}$\,.
 However,
 $$\barr{l}\giga\sum_{\chi'}\chi'(z\me)\fra1l\log\Big(\Psi(\fra{\lks(\chisa)}{\lki(\chi'\circ\ver)^l})\cdot\fra{\lks(\chi')^l}{\lks(\chi')^l}\Big)\\
 =l\cdot\giga\sum_{\chi'}\chi'(z\me)\fra1l\log\fra{\lks(\chi')}{\Psi\lki(\chi'\circ\ver)}+\giga\sum_{\chi'}
 \chi'(z\me)\fra1l\log\fra{\Psi\lks(\chisa)}{\lks(\chi')^l}
 \,,\earr$$ and, since $\psi_l\chi'={\chi'}^l$ here, the second
 summand is
 $$\barr{l}\giga\sum_{\chi'}\chi'(z\me)\fra1l\log
 \Big(\fra{\Psi\lks(\chisa)}{\lks(\chi')^l}\cdot\fra{\Psi\lks(\plcs)}{\Psi\lks({\chi'}^l)}\Big)\\
 =\giga\sum_{\chi'}\chi'(z\me)\fra1l\log\fra{\Psi\lks(\plcs)}{\lks(\chi')^l}+0=-\fra1l
 \ol{(t_{\Ki/k'}|\Ga)(z\me)}\ol z\,,\earr$$ by 2.~of Lemma 4.3 and
 Proposition 2.3.
 \mn Putting things together, we arrive at
 $$\barr{l}\ol{(\tki|\Ga)(z\me)}\ol
 z=\fra1l\ol{(t_{\Ki/k'}|\Ga)(z\me)}\ol z
 +\fra1{[G':\Ga]}\sum_{\chi'}\chi'(z\me)\fra1l\log\fra
 {\lks(\chi')}{\Psi\lki(\chi'\circ\ver)}\ -\\ \quad-\fra1l\ol{(t_{\Ki/k'}|\Ga)(z\me)}\ol
 z
 =\fra1{[G':\Ga]}\sum_{\chi'}(z\me)\fra1l\log\fra{\lks(\chi')}{\Psi\lki(\chi'\circ\ver)}\
 .\earr$$ Now insert condition 3.~and the above corollary with
 $x=\fra{\la_{\Ki/k'}}{\ver(\la_{\kab/k})}-1$\,; recall that $\Det(1+x)(\chi')=(\Det\fra{\la_{\Ki/k'}}{\ver(\la_{\kab/k})})(\chi')=
 \fra{\lks(\chi')}{\Psi\lki(\chi'\circ\ver)}$\,.
 This finishes the proof of the proposition.
 \Section{5}{Appendix}
 We begin this section by requiring only that $S$ contains the
 infinite primes and those above\nolinebreak[4] $l$.
 \Proposition{5.1} The existence of a pseudomeasure in $K_1(\qgi)$ is
 independent of the choice of the set $S$ as above. \Stop
 The proof of the proposition starts from the definition of a
 ``$K_1$-Euler factor for $\fp$ in $K/k$'' where $\fp$ is a prime
 of $k$ which does not divide $l\infty$. In this case, if $\fP$ is a
 fixed prime of $K$ above $\fp$, its decomposition group $D$ is an
 open subgroup of $G$ and its ramification subgroup $I$ a finite
 normal subgroup in $D$. We define
 \ben\item if $\fp$ is undecomposed, i.e., $D=G$,
 $E(\fp,K/k)\df[\qgi,1-\fra{g}{N\fp}\ve]\in K_1(\qgi)$\,, where
 $\ve=\fra1{|I|}\sum_{h\in I}h\,,\,N\fp=|\fo_k/\fp|$ and $gI$ is
 the Frobenius automorphism $\frp$ for the unique prime $\fP$ of
 $K$ above $\fp$\,, \item in general,
 $\eulerk\df\ind_D^G(\euler1)\in K_1(\qgi)$ with $\fp_1=\fP\cap
 k_1$ and $k_1=K^D$.\een
 Since the primes of $K$ above $\fp$ are $G$-conjugates of $\fP$,
 this definition does not depend on the choice of $\fP$.
 \Lemma{5.2} Let $G'$ be an open subgroup of $G$ with fixed field
 $k'$. Then $\res_G^{G'}\eulerk=\prod_{\fp'|\fp}\eulers$ with the $\fp'$
 running through all primes of $k'$ above $\fp$.\Stop
 \proof (of the lemma). We first look at the case when $D=G$ and
 use coset representatives $x_j\ (1\le j\le e\,,\,x_1=1)$ of
 $I\cap G'$ in $I$, so $I=\bigcup_{j=1}^e x_j(I\cap G')$, for
 getting the $\qgs$-basis $\{g^i\ve,(x_j-1)g^i:0\le i<f\,,\,1<j\le
 e\}$ of $\qgi$. Here, $e,f$ are the ramification index, respectively
 residue degree, of $\fp$ in $k'/k$.
 \sn Right multiplication by $1-\fra{g}{N\fp}\ve$ then yields the
 $e\times e$ block matrix, with blocks of size $f\times f$, $$\bmx{cccc}M&&&\\ &1&&\\
 &&\ddots&\\ &&&1\emx\quad\mr{where}\quad
 M=\bmx{cccc}1&-\fra1{N\fp}&&\\
 &\ddots&\ddots&\\&&1&-\fra1{N\fp}\\-\fra{g'}{N\fp}&&&1\emx$$ and
 $g^f=g'h$ for some $g'\in G'$ and $h\in I$. Of course, $g'\cdot
 (I\cap G')$ is the Frobenius automorphism $\frs$ of $\fP$ in $K/k'$. It
 follows that $\res_G^{G'}[\qgi,1-\fra{g}{N\fp}\ve]\in K_1(\qgs)$
 is represented by the above block matrix $M$ and thus equals (compare
 [CR, \S40C])
 $[\qgs,1-\fra{g'}{(N\fp)^f}]=[\qgs,1-\fra{g'}{N\fp'}]$\,, as is
 seen by performing the obvious column operations on $M$ to
 arrive at $$\bmx{ccccc}1&&&&\\ &1&&&\\&&&\ddots&\\
 \ast&\dots&&\ast&1-\fra{g'}{(N\fp)^f}\emx\,.$$
 \bn The general case of the lemma is reduced to this by means of
 the Mackey formalism:
 \bqo The double coset decomposition $G=\bigcup_{t\in T}G'tD$
 yields $$\res_G^{G'}\ind_D^G=\prod_{t\in T}\ind_{G'\cap
 D_t}^{G'}\res_{D_t}^{G'\cap D_t}c_t$$ where $D_t=tDt\me$ and
 $c_t:K_1({\cal{Q}}D)\to K_1({\cal{Q}}D_t)$ is induced by
 $x\mapsto txt\me$.\eqo
 We assume that $1\in T$ and start from the fixed prime $\fP$ in
 $K$ above $\fp$. Then $t\fP\cap k'$, for $t\in T$, are precisely
 the primes of $k'$ above $\fp$. Denote the fixed field of $D_t$
 by $k_t$ and set $\fp_t=t\fP\cap k_t$, so $\fp_t$ is undecomposed
 in $K/k_t$. Then
 $$\barr{l}\res_G^{G'}\eulerk=\res_G^{G'}\ind_D^G\euler1=\prod_t\ind_{G'\cap
 D_t}^{G'}\res_{D_t}^{G'\cap D_t}c_t\euler1\\
 =\prod_t\ind_{G'\cap D_t}^{G'}\res_{D_t}^{G'\cap
 D_t}\eulert=\prod_t\ind_{G'\cap D_t}^{G'}\eulerts\earr$$ with
 $k'_t$ the fixed field of $G'\cap D_t$ and $\fp'_t=t\fP\cap
 k'_t$\,, by the first part of the proof. Since $\fp'_t\cap k'=t\fP\cap k'$ divides $\fp$ and since
 the decomposition group of $t\fP$ in $K/k'$ is $G'\cap D_t$, the
 above product is $\prod_t E(t\fP\cap k',K/k')$, by the definition of the
 $K_1$-Euler factor in the general case.
 \mn This finishes the proof of the lemma and we turn back to the
 proof of Proposition 5.1 for which it now suffices to show that
 $\fp\notin S$ implies $$\Det\,\eulerk=\frac{L_{K/k,S\cup\{\fp\}}}{L_{K/k,S}}\,.$$
 We check this by evaluating both sides of the claimed equality at
 the characters $\chi$ of $G$. By Brauer induction we
 may assume that $\chi=\ind_{G'}^G\chi'$ is induced from a degree
 1 character of some open subgroup $G'$ of $G$. Now, with an
 obvious notation, Lemma 5.2 gives
 $$\barr{l}(\Det\,\eulerk)(\chi)=(\Det(\res_G^{G'}\eulerk))(\chi')=
 \prod_{\fp'|\fp}(\Det\,\eulers)(\chi')\\
 =
 \prod_{\fp'|\fp}\Det[\qgs,1-\fra{g_{\fp'}}{N\fp'}\ve_{\fp'}](\chi')\,.\earr$$
 Since $1-\fra{g_{\fp'}}{N\fp'}\ve_{\fp'}=1-\fra1{N\fp'}\cdot\fra1{|I_{\fp'}|}\sum_{h\in
 I_{\fp'}}g_{\fp'}h$\,, the value of
 $\Det[\qgs,1-\fra{g_{\fp'}}{N\fp'}\ve_{\fp'}]$ at the
 1-dimensional $\chi'$ is
 $=1-\fra1{N\fp'}\chi'(g_{\fp'})\chi'(\ve_{\fp'})\ol g_{\fp'}$
 (compare the first two displayed formulas following {\sc Fact 1}
 in \S1 and note that in the 1-dimensional case Det and Tr
 coincide). The above expression equals
 $1-\fra{\chi'(g_{\fp'})}{N\fp'}\ol g_{\fp'}$ or 1 according as
 $I_{\fp'}\subset\ker\chi'$ or not. Thus
 $$(\Det\,\eulerk)(\chi)=\prod_{\fp'|\fp,\,I_{\fp'}\subset\ker\chi'}(1-\fra{\chi'(\frs)}{N\fp'}\ol\frs)$$
 and we are left with showing that this product is also
 $$=\ \frac{L_{K/k,S\cup\{\fp\}}(\chi)}{L_{K/k,S}(\chi)}=\fra{L_{K/k',S'\cup\{\fp'\}}(\chi')}{L_{K/k',S'}(\chi')}$$
 with $\{\fp'\}$ again the set of primes of $k'$ above $\fp$.
 \sn Now, with $u'\in 1+l\zl$ denoting the action of a generator
 $\ga_{k'}$ of $\Ga_{k'}$ on $l$-power roots of unity, the ratio
 of power series
 $\fra{G_{\chi',S'\cup\{\fp'\}}(T)}{G_{\chi',S'}(T)}$ is uniquely
 determined by its values
 $$\barr{l}\fra{G_{\chi',S'\cup\{\fp'\}}((u')^n-1)}{G_{\chi',S'}((u')^n-1)}=
 \fra{L_{l,S'\cup\{\fp'\}}(1-n,\chi')}{L_{l,S'}(1-n,\chi')}=\prod_{\fp'|\fp,\,I_{\fp'}\subset\ker\chi'}
 (1-\fra{\chi'(\frs)}{N\fp'}(N\fp')^n)\\
 =\prod(1-\fra{\chi'(\frs)}{N\fp'}\group{N\fp'}^n)=\prod(1-\fra{\chi'(\frs)}{N\fp'}(u')^{nb_{\fp'}})\earr$$
 at all natural numbers $n\equiv0\mod l-1$. Here,
 $\group{N\fp'}\in1+l\zl$ is determined by
 $(N\fp')^n=\group{N\fp'}^n$ for all such $n$ and $b_{\fp'}\in\zl$
 by $\group{N\fp'}=(u')^{b_\fp'}$.
  Consequently, substituting $T=\ga_{k'}-1$,
 $$\fra{L_{K/k',S'\cup\{\fp'\}}(\chi')}{L_{K/k',S'}(\chi')}=
 \fra{G_{\chi',S'\cup\{\fp'\}}(\ga_{k'}-1)}{G_{\chi',S'}(\ga_{k'}-1)}=\prod(1-\fra{\chi(\frs)}{N\fp'}
 \ga_{k'}^{b_{\fp'}})$$ and so $\ol\frs=\ga_{k'}^{b_{\fp'}}$
 finishes the proof of the proposition\,:
 these automorphisms act on $l$-power roots of unity by $\group{N\fp'}=(u')^{b_{\fp'}}$.
 \bn We close this section by recalling the ``equivariant main
 conjecture''. The set $S$ is now supposed to be sufficiently
 large. Moreover, we assume $\mu=0$ for $\Ki/k$
 \footnote{the assumption is independent of $S$ as long as $l\in S$, see [NSW, (11.3.6), p.615]}.
 \sn Let $M$ be the maximal abelian
 $l$-extension of $\Ki$, which is unramified outside $S$, and set
 $X=G_{M/\Ki}.$ Then $X$ is a finitely
 generated torsion $\lgi$-module \footnote{As before, ``torsion'' means that
 there is a central regular $c\neq0$ in $\lgi$ which annihilates
 $X$.} (the so-called Iwasawa module). Though $X$
 generally does not have finite projective dimension itself, it naturally induces
 an element \footnote{We prefer to just write $\mho$ and $\Theta$ rather than
 $\ti\mho,\ti\Theta$, as in [RW\,2,3,4]. Note however that these $\mho,\Theta$ differ
 from the ones appearing in [RWt,RW1]: compare [RW2, Remark A., p.564].}
 $\mho$ in the Grothendieck group $K_0T(\lgi)$ of all
 finitely generated torsion $\lgi$-modules with finite projective
 dimension. This $\mho$ not only keeps all the information of the
 Iwasawa module $X$, but also includes extension class
 data. For its derivation see [RWt or RW2].
 \mn The localization sequence of $K$-theory provides a connecting
 homomorphism $$K_1(\qgi)\sr{\partial}{\lto}K_0T(\lgi)$$ and $\partial$ has
 $\mho$ in its image [RW2, Lemma 13]. Combining $\partial$ and Det
 yields $$\barr{l} \pht{xxxx}K_1(\qgi)\sr\partial\lto K_0T(\lgi)\\
 \pht{xxxx}\daz{\Det}\\
 \Hom^\ast(\rlgi,(\qcgak)\mal)\earr$$ and the question arises whether there is a
 common source in $K_1(\qgi)$ for the two distinguished elements,
 $\mho$ at the right and $\lki$ at the bottom. The ``equivariant
 main conjecture'' of Iwasawa theory, as stated in [RW2], asserts
 that, roughly speaking, there is a preimage $\Theta$ of $\lki$ in
 $K_1(\qgi)$ such that $\partial(\Theta)=\mho$\,. If $\gi$ is
 abelian, this is essentially what the Main Conjecture of
 classical Iwasawa theory is about (compare [RW1]). In [RW3,4] it
 has been shown that, if $\mu=0$, the existence of $\Theta$ is
 equivalent to $\lki$ belonging to $\Det\,K_1(\lbgi)$. Then
 $\Theta$ is a non-abelian pseudomeasure for $\Ki/k$.
 \enlargethispage{5mm}
 \Proposition{5.3} The ``equivariant main conjecture'' does not
 depend on the choice of a sufficiently large set $S$.\Stop
 It is enough to show
 \nf{$\ast$}{\partial\,\eulerk=\mho_{S\cup\{\fp\}}-\mho_S}
 for $\fp\notin S$. In the notation of the general case of the
 definition of a $K_1$-Euler factor we have
 $$\eulerk=\ind_D^G\euler1=\ind_D^G[{\cal{Q}}D,1-\fra{g_1}{N\fp_1}\ve_1]=\ind_D^G[{\cal{Q}}D\ve_1,1-\fra{g_1}{N\fp_1}]$$
 because $1-\fra{g_1}{N\fp_1}\ve_1$ acts on
 ${\cal{Q}}D=({\cal{Q}}D)\ve_1\oplus({\cal{Q}}D)(1-\ve_1)$ as
 $(1-\fra{g_1}{N\fp_1})\ve_1+(1-\ve_1)$. Since
 $({\cal{Q}}D)\ve_1={\cal{Q}}(D/I)$, $N\fp_1=N\fp$, and as $g_1$ acts
 here as $\mr{Fr}_{\fp_1}$ this becomes
 $$\partial\,\eulerk=\partial\Big(\ind_D^G[{\cal{Q}}(D/I),1-\fra{\mr{Fr}_{\fp_1}}{N\fp}]\Big)\,.$$
 To compare this with the right hand side of $(\ast)$ we now apply
 the proof of [RWt, Proposition 4.7] (which is not built on Leopoldt's conjecture). Because $S$ is sufficiently large, this only requires us to
 restate the second and fourth displayed formula on [loc.\,cit.,
 p.38] as
 $$\mho_{S\cup\{\fp\}}-\mho_S=[\ind_D^GC]=\partial\Big(\ind_D^G[{\cal{Q}}(D/I),1-\fra{\mr{Fr}_{\fp_1}}{N\fp}]\Big)\,.$$
 Note that in [RWt] we inverted only the ``never zero divisors''
 $R\subset\La G$, so these assertions remain true on inverting the
 larger set of regular central elements.
  \bbn{\large {\sc References}}
 \small
 \bbn
 \btb{rp{13cm}}
 \,[CR]   & Curtis, C.W.~and Reiner, I., {\em Methods of
            Representation Theory,
            I,II.} John Wiley \& Sons (1981,1987) \\
 \,[C-N]  & Cassou-Nogu\`es, P., {\em Valeurs aux entiers n\'egatifs des fonctions zêta et fonctions zêta
            $p$-adiques.} Invent.~Math. {\bf 51} (1979), 29-59\\
 \etb
 \noi
 \btb{rp{13cm}}
 \,[DR]   & Deligne, P.~and Ribet, K., {\em Values of abelian
            $L$-functions at negative integers
            over totally real fields.} Invent.~Math. {\bf 59} (1980), 227-286\\
 \,[Fr]   & Fröhlich, A., {\em Galois Module Structure of Algebraic
            Integers.} Springer-Verlag (1983)\\
 \,[Gr]   & Greenberg, R., {\em On $p$-adic Artin $L$-functions.} Nagoya Math.~J. {\bf 89} (1983), 77-87\\
 \,[NSW]  & Neukirch, J., Schmidt, A.~and Wingberg, K., {\em
            Cohomology of Number Fields.} Springer Grundlehren der math.~Wiss.~{\bf 323} (2000)\\
  \,[RWt]  & Ritter, J.~and Weiss, A., {\em The Lifted Root Number
            Conjecture and Iwasawa Theory.} Memoirs of the AMS {\bf 748} (2002)\\
 \,[RW1]  & ------------------------\,, {\em Toward equivariant Iwasawa theory.}
            Manuscripta Math. {\bf 109} (2002), 131-146\\
 \,[RW2,3,4]  & ------------------------\,, {\em Toward equivariant Iwasawa theory, II; III; IV.}
           Indagationes Mathematicae {\bf 15} (2004), 549-572; Mathematische Annalen {\bf 336}
           (2006), 27-49 [\,DOI:
            10.1007/s00208-006-0773-4\,]; Homology, Homotopy and Applications {\bf 7} (2005), 155-171
            [\,http://intlpress.com/HHA/v7/n3/\,] \\
 \,[Se]   & Serre, J.-P., {\em Sur le r\'esidu de la fonction
           z\^{e}ta $p$-adique d'un corps de nombres.}
           C.R.Acad.Sci.~Paris {\bf 287} (1978), s\'erie A, 183-188\\
 \,[Ty1]  & Taylor, M.J., {\em Galois module type congruences for values of $L$-functions.} J.~LMS {\bf
            24} (1981), 441-448\\
 \,[Ty2]  & ------------, {\em Classgroups of Group Rings.} LMS
            Lecture Notes Series {\bf 91}, Cambridge (1984)\\
 \etb

 \vsp{.5} {\footnotesize \bct Institut für Mathematik $\cdot$
 Universität Augsburg $\cdot$ 86135 Augsburg $\cdot$ Germany \\
 Department of Mathematics $\cdot$ University of Alberta $\cdot$
 Edmonton, AB $\cdot$ Canada T6G 2G1   \ect

 }

 \end{document}